\documentclass[a4paper, 11pt]{amsart}
\usepackage[utf8]{inputenc}
\usepackage[T1]{fontenc}
\usepackage[english]{babel}
\usepackage[margin=1in]{geometry} % fullpage
\usepackage{microtype, xcolor, graphicx}
\usepackage{amsmath, amssymb, amsthm}
\usepackage[linktocpage=true]{hyperref}
\usepackage[capitalise]{cleveref}
\usepackage{enumitem}
\usepackage{tikz-cd} % https://tikzcd.yichuanshen.de

\theoremstyle{plain}
\newtheorem{mainthm}{Theorem}

\newtheorem{thm}{Theorem}[section]
\newtheorem{lem}[thm]{Lemma}
\newtheorem{prop}[thm]{Proposition}
\newtheorem{cor}[thm]{Corollary}

\theoremstyle{definition}
\newtheorem{defn}[thm]{Definition}
\newtheorem{exmp}[thm]{Example}
\newtheorem{qst}[thm]{Question}

\theoremstyle{remark}
\newtheorem{rem}[thm]{Remark}

\numberwithin{equation}{section}

\DeclareMathOperator{\Pos}{Pos}
\DeclareMathOperator{\Mods}{Mod}
\DeclareMathOperator{\mods}{mod}
\DeclareMathOperator{\Hom}{Hom}

\DeclareMathOperator{\otimesL}{\otimes_R^\mathbf{L}}
\DeclareMathOperator{\chom}{\raisebox{0.5pt}{$\mathcal{H}$}om}

\DeclareMathOperator{\Spec}{Spec}
\DeclareMathOperator{\Spc}{Spc}
\DeclareMathOperator{\Supp}{Supp}

\DeclareMathOperator{\id}{id}
\DeclareMathOperator{\coker}{coker}

\DeclareMathOperator{\hocolim}{\underrightarrow{\mathrm{hocolim}}}
\DeclareMathOperator{\proj}{proj}

\DeclareMathOperator{\thick}{\mathsf{thick}}
\DeclareMathOperator{\loc}{\mathsf{loc}}
\DeclareMathOperator{\susp}{\mathsf{susp}}

\DeclareMathOperator{\Filt}{Filt}
\DeclareMathOperator{\Serre}{Serre}
\DeclareMathOperator{\Wide}{Wide}
\DeclareMathOperator{\Thick}{Thick}
\DeclareMathOperator{\Loc}{Loc}
\DeclareMathOperator{\Susp}{Susp}
\DeclareMathOperator{\Aisle}{Aisle}
\newcommand{\cA}{\mathcal{A}}

\newcommand{\cC}{\mathcal{C}}
\newcommand{\cD}{\mathcal{D}}

\newcommand{\cK}{\mathcal{K}}
\newcommand{\cL}{\mathcal{L}}
\newcommand{\cM}{\mathcal{M}}
\newcommand{\cN}{\mathcal{N}}
\newcommand{\cO}{\mathcal{O}}
\newcommand{\cP}{\mathcal{P}}

\newcommand{\cS}{\mathcal{S}}
\newcommand{\cT}{\mathcal{T}}
\newcommand{\cU}{\mathcal{U}}
\newcommand{\cV}{\mathcal{V}}

\newcommand{\cX}{\mathcal{X}}

\newcommand{\cZ}{\mathcal{Z}}

\newcommand{\m}{\mathfrak{m}}
\newcommand{\p}{\mathfrak{p}}
\newcommand{\q}{\mathfrak{q}}

\newcommand{\bZ}{\mathbb{Z}}

\newcommand{\bK}{\mathbb{K}} % Used for fields

\newcommand{\longmapsfrom}{\mathrel{\reflectbox{\ensuremath{\longmapsto}}}}
\newcommand{\leftangle}{\left\langle}
\newcommand{\rightangle}{\right\rangle}
\newcommand{\leftcurly}{\left\{}
\newcommand{\rightcurly}{\right\}}

\begin{document}
\title[Balmer spectrum and tensor telescope conjecture]{The Balmer spectrum and tensor telescope conjecture for noetherian path algebras}
\author{Enrico Sabatini}
\address{Enrico Sabatini, Dipartimento di Matematica ``Tullio Levi-Civita'', Universit{\`a} degli Studi di Padova, via Trieste 63, 35121 Padova, Italy}
\email{enrico.sabatini@phd.unipd.it}
\thanks{The author was supported by the Department of Mathematics ``Tullio Levi-Civita'' of the University of Padova}
\begin{abstract}
    Given a commutative noetherian ring $R$ and a finite acyclic quiver $Q$, we study the tensor triangulated category $\cD(RQ)$ endowed with the vertexwise tensor product. We find a description of the internal hom functor and show that the category is not rigid. We compute its Balmer spectrum and, despite the non-rigidity, we get a classification of all the thick tensor-ideals and a stratification result. After introducing the notion of tensor-t-structure, we give a classification of the compactly generated ones and prove the tensor telescope conjecture.
\end{abstract}
\maketitle
\setcounter{tocdepth}{1}
\tableofcontents

\section*{Introduction}
Tensor triangulated categories (\emph{tt-categories}) arise naturally in homotopy theory, algebraic geometry and modular representation theory. They appear, for example, as ``unital algebraic stable homotopy categories'' in \cite{HPS}, ``symmetric monoidal presentable stable $\infty$-categories'' in \cite{Lur} and ``stable module categories of finite groups'' in \cite{BIK}. A seminal contribution to this field is due to Balmer \cite{Bal} and Balmer-Favi \cite{BF}, where ``small'' and ``big'' tt-categories are considered -- namely, \emph{essentially small} tt-categories and \emph{rigidly-compactly generated} tt-categories. In particular, every ``big'' tt-category yields a ``small'' one via its full triangulated subcategory of (rigid-)compact objects. The most interesting feature of a ``small'' tt-category is the \emph{Balmer spectrum}, which, analogously to the prime spectrum of a commutative ring, is both a topological space and a poset encoding structural information about the category. It provides a notion of support for objects and allows for the classification of relevant tensor-compatible subcategories. As for ``big'' tt-categories, thanks to \emph{rigidity}, the Balmer spectrum of the subcategory of compact objects provides a well-founded notion of support also for objects that are not necessarily compact. The ability of the Balmer spectrum to classify the relevant tensor-compatible subcategories of a ``big'' tt-category is known as \emph{stratification}, and it is not always satisfied (see \cite{BHS} for a foundational work on this topic).

In this paper we consider the derived category of representations of a finite quiver over a commutative noetherian ring endowed with the (derived) vertexwise tensor product. This turns out to be a tensor triangulated category. In this framework, we compute the \emph{Balmer spectrum} of its compact objects and derive several \emph{tt-classifications}. This computation is noteworthy because, although the category is compactly generated, its compact objects do not form a rigid tt-category. Explicit computations of Balmer spectra in such non-rigid situations are rare -- see \cite{LS, Xu, MT, BCS} -- and, among these, only the first two consider compact objects of a non-rigidly compactly generated tt-category.

The relevance of our example also lies in the fact that, despite the lack of rigidity, it is still the subcategory of compact objects that captures the relevant information about the ambient category. Usually, this is the expected behavior in rigidly-compactly generated tt-categories $\cT$, where the subcategory of compact objects $\cT^c$ coincide with the subcategory of rigid objects $\cT^d$. However we are not in such situation. In fact, when $\cT^c$ fails to be rigid -- or even to form a tt-subcategory, as it may happen when $\cT$ is not a compactly generated tt-category -- it is generally regarded as more appropriate, as suggested by \cite{BCHS}, to consider the spectrum of rigid objects $\Spc(\cT^d)$ rather than the one of the compacts. For example, \cite{Ken} exhibits big tt-categories where either $\cT^c\subsetneq\cT^d$ or $\cT^d\subsetneq\cT^c$, and in both instances the meaningful spectrum is $\Spc(\cT^d)$. Our example suggests that the recourse to rigid objects may be relevant only when $\cT$ is not a compactly generated tt-category, e.g.~when $\cT^c\neq\cT^d$.

We show that the Balmer spectrum $\Spc(\cD^c(RQ))$ is a disjoint union of copies of the prime spectrum $\Spec(R)$, extending the result of \cite{LS}, and that its specialization-closed subsets classify all thick tensor-ideals of $\cD^c(RQ)$. 

\begin{mainthm}[\cref{Spectrum} (2) and \cref{Thickbyspc}]\label{ThmA}
    For any commutative noetherian ring $R$ and finite acyclic quiver $Q=(Q_0,Q_1)$, there is an order-reversing homeomorphism
    \[f:\Spc(\cD^c(RQ))\longrightarrow\Spec(R)\times Q_0\]
    i.e.~$\Spc(\cD^c(RQ))\overset{f}{\cong}_\mathrm{\,Top}\Spec(R)\times Q_0$ and $\Spc(\cD^c(RQ))\overset{f}{\cong}_{\,\Pos}\Spec(R)^\mathrm{op}\times Q_0$, where $Q_0$ is taken with discrete topology and trivial ordering.
    
    Moreover, thick tensor-ideals of $\cD^c(RQ)$ are in bijection with specialization closed subsets of the Balmer spectrum $\Spc(\cD^c(RQ))$.
\end{mainthm}

Note that this computation is also significant from the viewpoint of abstract tensor triangular geometry. Indeed, in the rigid case, if the spectrum decomposes as a disjoint union of subspaces, then the category itself splits as a direct sum of subcategories indexed by those components (see \cite[Theorem 5.14]{BF}), but our category is indecomposable.
\smallskip

Moving from the compact objects to the ambient tt-category, we explain how the absence of rigidity obstructs the construction of the Balmer-Favi support, and we introduce an \emph{ad hoc} support (see \cref{BigSupp}) which we prove to be \emph{stratifying}. 

\begin{mainthm}[\cref{Strat} (2)]\label{ThmB}
    For any commutative noetherian ring $R$ and finite acyclic quiver $Q$, taking the support of objects in $\cD(RQ)$ defines a bijection between localizing tensor-ideals of $\cD(RQ)$ and subsets of the Balmer spectrum $\Spc(\cD^c(RQ))$
    \[\begin{tikzcd}
        \Loc_\boxtimes(\cD(RQ)) \arrow[rr, "\Supp_{RQ}", shift left] & & \leftcurly\text{Subsets of }\Spc(\cD^c(RQ))\rightcurly \arrow[ll, "\Supp_{RQ}^{-1}", shift left]
    \end{tikzcd}\]
\end{mainthm}

Next, we turn to the recently introduced notion of tensor-t-structures, defined with respect to a fixed tensor-closed suspended subcategory $\cT^{\leq 0}$ containing the tensor unit. In this framework, one can formulate the \emph{tensor telescope conjecture}, which asks whether every homotopically smashing tensor-t-structure is compactly generated. We provide a classification of compactly generated tensor-t-structures with respect to the standard aisle $\cD^{\leq 0}(RQ)$ in terms of filtrations of specialization-closed subsets of the Balmer spectrum (see \cref{Aislebyspc}), and we establish the following result:

\begin{mainthm}[Theorems \ref{TTCthm} and \ref{TTCimprove}]\label{ThmC}
    For any commutative noetherian ring $R$ and finite acyclic quiver $Q$, the derived category $\cD(RQ)$ satisfies the telescope conjecture for tensor-t-structures with respect to $\cD^{\leq0}(RQ)$. More generally, the same holds for tensor-t-structures with respect to any tensor-closed suspended subcategory generated by a filtration system with Dynkin support (see \cref{FiltSyst}).
\end{mainthm}

This result not only generalizes \cite[Theorem 4.2]{Sab} -- which plays a crucial role in its proof -- but also marks a significant step toward establishing the telescope conjecture for representations of finite acyclic simply-laced quivers over commutative noetherian rings.

\subsection*{Structure of the paper}
The paper is divided into four sections. In \textbf{\cref{Sec1}}, we introduce tensor triangulated categories, provide an overview of the tensor triangulated structure on the derived category $\cD(RQ)$, and conclude with some properties of derivators that will be extensively used in the sequel. In \textbf{\cref{Sec2}}, we review the basics of the Balmer spectrum, analyze it in our context, and prove Theorems \ref{ThmA} and \ref{ThmB}. In \textbf{\cref{Sec3}}, we define tensor t-structures and characterize them in our setting. Finally, in \textbf{\cref{Sec4}}, we establish the tensor telescope conjecture of \cref{ThmC}.

\subsection*{Acknowledgments}
The author would like to thank his supervisors, Jorge Vit{\'o}ria and Michal Hrbek: the former for carefully reading a preliminary version of this paper and providing valuable feedback, and the latter for his hospitality during the author's visit to the Czech Academy of Sciences and for engaging in long and helpful discussions on the topics of this work. In particular, the author thanks Michal for suggesting the topic addressed in \cref{StratSec} and for pointing out the arguments used in the proof of \cref{SetGen}.

\section{Preliminaries}\label{Sec1}
\subsection{Tensor triangulated categories}
Let $(\cT,\otimes,\mathbf{1})$ be a \emph{tensor triangulated category} (\emph{tt-category}) in the sense of \cite{Bal}, i.e.~a triangulated category $\cT$ with a symmetric monoidal product $\otimes:\cT\times\cT\to\cT$ with unit $\mathbf{1}\in\cT$ and such that for any $X\in\cT$ the functor $X\otimes\_:\cT\to\cT$ is triangulated. Recall that $\cT$ is a \emph{compactly generated tt-category} if:
\begin{itemize}
    \item[(i)] $\cT$ is a compactly generated triangulated category, i.e.~it has coproducts and the subcategory $\cT^c$ of its compact objects is skeletally small and generates $\cT$ as a localizing subcategory;
    \item[(ii)] The tensor product $\otimes:\cT\times\cT\to\cT$ preserves coproducts;
    \item[(iii)] $\cT^c$ is a \emph{tt-subcategory}, i.e.~the tensor product restricts to $\otimes:\cT^c\times\cT^c\to\cT^c$ and $\mathbf{1}\in\cT^c$.
\end{itemize}
In this setting, by Brown representability \cite[Theorem 4.1]{Nee}, the tensor product $\otimes:\cT\times\cT\to\cT$ is a \emph{closed} monoidal product, that is the functor $X\otimes\_$ admits a right adjoint $[X,\_]:\cT\to\cT$ for all $X\in\cT$. We call the bifunctor $[\_,\_]:\cT\times\cT\to\cT$ the \emph{internal hom} of $\cT$. Following \cite{BF}, a compactly generated tt-category $\cT$ is called \emph{rigidly-compactly generated} if:
\begin{itemize}
    \item[(iv)] The internal hom restricts to $[\_,\_]:\cT^c\times\cT^c\to\cT^c$;
    \item[(v)] $\cT^c$ is a \emph{rigid tt-category}, i.e.~for any $X,Y\in\cT^c$ the \emph{evaluation map} $[X,\mathbf{1}]\otimes Y\to[X,Y]$ is an isomorphism.
\end{itemize}
The objects $X\in\cT$ satisfying the property (v) for any $Y\in\cT$ are called \emph{rigid} objects of $\cT$. Note that, by \cite[Theorem 2.1.3 (d)]{HPS}, the condition (v) implies that compact objects and rigid objects coincide in $\cT$. Moreover, denoting by $X^\vee$ the object $[X,\mathbf{1}]$, by \cite[III.1.3 (i)]{LMS}, any rigid object $X$ is \emph{dualizable}, i.e.~$X\cong(X^\vee)^\vee$, and it is a direct summand of $X^\vee\otimes X\otimes X$.
\medskip

From now on, let $R$ be a commutative noetherian ring, $Q$ a finite acyclic quiver and denote by $RQ$ the path algebra of $Q$ over $R$. Recall that the category of modules $\Mods(RQ)$ is equivalent to the category of representations of $Q=(Q_0,Q_1)$ by $R$-modules, i.e.~the category with objects $X=(X_i,X_\alpha)_{i\in Q_0,\,\alpha\in Q_1}$, where $X_i$ is an $R$-module for any vertex $i\in Q_0$ and $X_\alpha:X_i\to X_j$ is an $R$-linear map for any arrow $\alpha:i\to j\in Q_1$, and morphisms $f:X\to Y$ where $f=(f_i:X_i\to Y_i)_{i\in Q_0}$ such that $f_j\circ X_\alpha=Y_\alpha\circ f_i$ for any arrow $\alpha:i\to j\in Q_1$. 

Given two $RQ$-modules $X=(X_i,X_\alpha)$ and $Y=(Y_i,Y_\alpha)$, we define a symmetric tensor product as $X\boxtimes_{RQ}Y=(X_i\otimes_R Y_i,X_\alpha\otimes_R Y_\alpha)$. In the following we will consider the tt-category $(\cD(RQ),\boxtimes_{RQ}^\mathbf{L},\mathbf{U})$, where $\cD(RQ)$ is the derived category of the module category $\Mods(RQ)$, the symmetric monoidal structure is given by the left derived functor of $\boxtimes_{RQ}$ and the tensor unit $\mathbf{U}(=\mathbf{U}[0])$ is the stalk complex defined by $\mathbf{U}_i=R$ and $\mathbf{U}_\alpha=\id_R$ for any $i\in Q_0$ and $\alpha\in Q_1$.

\begin{rem}
    Note that the tensor unit is not a generator of the triangulated category, indeed $\loc_{RQ}\langle\mathbf{U}\rangle\neq\cD(RQ)$. There are many results available for tt-categories $(\cT,\otimes,\mathbf{1})$ such that $\loc_\cT\langle\mathbf{1}\rangle=\cT$, which we can not therefore use in our context. In particular, a localizing (resp.~thick) subcategory of $\cD(RQ)$ (resp.~$\cD^c(RQ)$) will not be automatically a localizing (resp.~thick) tensor-ideal of $\cD(RQ)$ (resp.~$\cD^c(RQ)$) -- see \cref{Sec2} for definitions.
\end{rem}

\begin{prop}\label{CompGen}
    The tt-category $(\cD(RQ),\boxtimes_{RQ}^\mathbf{L},\mathbf{U})$ is compactly generated. In particular, $(\cD^c(RQ),\boxtimes_{RQ},\mathbf{U})$ is an essentially small tt-category.
\begin{proof}
    The derived category $\cD(RQ)$ is a compactly generated triangulated category with compact objects $\cK^b(\proj(RQ))$ and the tensor product $\boxtimes_{RQ}$ clearly commutes with coproducts. Recall from \cite[Theorem 3.1]{EE} that an $RQ$-module $P$ is projective if and only if, for any $i\in Q_0$, $P_i$ is projective and the map
    \[\bigoplus\limits_{\substack{\alpha\in Q_1\\t(\alpha)=i}}P_{s(\alpha)}\to P_i\]
    is a split monomorphism -- where $s(\alpha)$ and $t(\alpha)$ denote respectively the source and the target of the arrow $\alpha$. It follows easily that, given two projective $RQ$-module $P$ and $P'$, the tensor product $P\boxtimes_{RQ} P'$ is again projective. Since for any two compact objects $X,X'\in\cD^c(RQ)$, the product $X\boxtimes_{RQ}^\mathbf{L}X'$ is isomorphic to the $\oplus$-totalization of the double complex $X^i\boxtimes_R X'^j$, it is compact too. Finally, by \cite[Lemma 3.1 (i)]{CB}, the $RQ$-module $\mathbf{U}$ has projective dimension less or equal then $1$ and, being finitely presented, it is compact.
\end{proof}
\end{prop}

By Brown representability, it follows from \cref{CompGen} that the category $\cD(RQ)$ is endowed with an internal hom, which we will denote by $\chom_{RQ}(\_,\_):\cD(RQ)\times\cD(RQ)\to\cD(RQ)$.

For any $i\in Q_0$, let $P(i)$ be the projective $RQ$-module such that, for any vertex $k$, $P(i)_k$ is the free $R$-module with basis the set of all paths in $Q$ from $i$ to $k$ and, for any arrow $\beta:k\to\ell$, $P(i)_\beta:P(i)_k\to P(i)_\ell$ is given by left multiplication by $\beta$; for any $\alpha:i\to j\in Q_1$, let $\iota_\alpha:P(j)\to P(i)$ be the monomorphism given by right multiplication by $\alpha$. Then:

\begin{thm}
    Given left $RQ$-modules $X$ and $Y$, the internal hom is defined by
    \[\chom_{RQ}(X,Y):=\left(\Hom_{RQ}(X\boxtimes_{RQ}P(i),Y),\,\_\circ(\id_X\boxtimes_{RQ}\,\iota_\alpha)\right)_{i\in Q_0,\,\alpha\in Q_1}\]
\begin{proof}
    By \cite[Proposition IV.10.1]{Ste}, for any left $RQ$-module $X$, the functor $X\boxtimes_{RQ}\_$ is isomorphic to the functor $(X\boxtimes_{RQ}RQ)\otimes_{RQ}\_$. Here the module $X\boxtimes_{RQ}RQ$ is seen as an $RQ$-bimodule where, denoting by $\rho_\omega$ the right multiplication in $RQ$ by a path $\omega$, the right $RQ$-module structure of $X\boxtimes_{RQ}RQ$ is given by $\cdot\omega:=\id_X\boxtimes_{RQ}\,\rho_\omega$. Thus, the functor $X\boxtimes_{RQ}\_$ is left adjoint to the functor $\Hom_{RQ}(X\boxtimes_{RQ}RQ,\_)$. For any left $RQ$-module $Y$, the left $RQ$-module structure on $\Hom_{RQ}(X\boxtimes_{RQ}RQ,Y)$ is defined by $\omega\cdot(h_i)_{i\in Q_0}=(h_i\circ(\id_X\boxtimes_{RQ}\,\rho_\omega))_{i\in Q_0}$. In particular, at any vertex $i\in Q_0$ it is represented by the $R$-module $\varepsilon_i\cdot\Hom_{RQ}(X\boxtimes_{RQ}RQ,Y)=\Hom_{RQ}(X\boxtimes_{RQ}RQ\cdot\varepsilon_i,Y)=\Hom_{RQ}(X\boxtimes_{RQ}P(i),Y)$.
\end{proof}
\end{thm}

\begin{prop}
    The internal hom restricts to $\cD^c(RQ)$.
\begin{proof}
    For any two finitely generated projective $RQ$-modules $P$ and $P'$, by \cref{CompGen} the $RQ$-module $P\boxtimes_{RQ}P(i)$ is finitely generated projective too, and, by \cite[Theorem A]{CB}, $\Hom_{RQ}(P\boxtimes_{RQ}P(i),P')$ is a finitely generated projective $R$-module. Thus, $\chom_{RQ}(P,P')$ is an $RQ$-lattice and, by \cite[Lemma 3.1 (i)]{CB}, it is compact in $\cD(RQ)$. Since, for any two compact objects $X,X'\in\cD^c(RQ)$, the complex $\chom_{RQ}(X,X')$ is isomorphic to the $\oplus$-totalization of the double complex $\chom_{RQ}(X^{-i},X'^j)$, it is compact.
\end{proof}
\end{prop}

We will prove in the next proposition that the tt-category $(\cD^c(RQ),\boxtimes_{RQ},\mathbf{U})$ is not rigid whenever the quiver $Q$ is connected, has at least two vertices, and admits a source -- which is always the case for finite acyclic quivers. This result generalizes \cite[Remark 2.21]{Ito} and \cite[Example 6.13]{San} to any commutative noetherian ring $R$.

\begin{prop}\label{notrigid}
    For a non-trivial connected finite acyclic quiver $Q$, the tt-category $\cD^c(RQ)$ is not rigid.
\begin{proof}
    Any finite acyclic quiver has a source, say $s\in Q_0$. Let $U(s)$ be the representation such that, $U(s)_s=R$ and $U(s)_k=0$ for any $k\neq s$, note that it is compact by \cite[Lemma 3.1 (i)]{CB}. Since $s$ is a source, the product $U(s)\boxtimes_{RQ}P(i)$ is equal to $U(s)$ if $i=s$ or $0$ otherwise. Thus, since the quiver is non-trivial and connected, $\chom_{RQ}(U(s),\mathbf{U})=0$, while $\chom_{RQ}(U(s),U(s))=U(s)$. It follows that $U(s)$ is not rigid, i.e.~$\chom_{RQ}(U(s),\mathbf{U})\boxtimes_{RQ}U(s)\ncong\chom_{RQ}(U(s),U(s))$.
\end{proof}
\end{prop}

There are many results available for rigid tt-categories that we cannot apply in our context; in particular, not every thick tensor-ideal of $\cD^c(RQ)$ is automatically radical. Nevertheless, even though $\cD^c(RQ)$ is not rigid, we will show in \cref{Rigid} that every thick tensor-ideal is indeed radical, so the Balmer spectrum still provides a classification of all thick tensor-ideals. For further distinctions from the rigid setting, see \cref{StratSec}.

\subsection{Derivators}\label{DerSec}
Before we end this section, we want to introduce some tools coming from the theory of derivators, which will be used later. Despite this theory being very general, we will stick to our context as much is possible. Our main reference is \cite{Gro} and we refer to it for any unexplained terminology. By \cite[Proposition 1.30, Example 4.2]{Gro}, the assignment
\[\mathbb{D}_{R}:\mathrm{Cat^{op}}\to\mathrm{CAT} \text{ such that } \mathbb{D}_{R}(Q):=\cD(RQ) \text{ for any finite quiver } Q\]
form a strong and stable derivator. In particular, by axiom (Der3) in \cite[Definition 1.5]{Gro}, for any vertex $i\in Q_0$ there is an adjoint triple $i_!\dashv i^\ast\dashv i_\ast$ as follows
\[\begin{tikzcd} \cD(RQ) \arrow[rr, "i^\ast"] & & \cD(R) \arrow[ll, "i_!"', bend right] \arrow[ll, "i_\ast", bend left] \end{tikzcd}\]
where, by \cite[Corollary 4.19]{Gro}, all the functors are triangulated and $i^\ast$ is the \emph{evaluation functor} at $i$, induced by precomposition with the injection $\{i\}\hookrightarrow Q$. Note that, since coproducts in $\cD(RQ)$ are computed vertexwise and bounded complexes of finitely generated projective $RQ$-modules are vertexwise bounded complexes of finitely generated projective $R$-modules, the functor $i^\ast$ preserves both coproducts and compact objects, thus, by \cite[Theorem 5.1]{Nee}, $i_!$ preserves compacts and $i_\ast$ preserves coproducts. Moreover, $i_!$, $i^\ast$ and, by \cite[Theorem 4.1]{Nee}, $i_\ast$ are all left adjoint functors, so they preserve directed homotopy colimits (see \cite[Proposition 2.4, Corollary 2.12]{Gro}).

\begin{rem}\label{GroRem}
    By axiom (Der4) in \cite[Definition 1.5]{Gro} (see also \cite[Lemma 11.1]{KN}), for any complex $M\in\cD(R)$ and vertices $i,j\in Q_0$ we have that
    \[j^\ast i_! M\cong\coprod_{Q(i,j)}M \text{ and } j^\ast i_\ast M\cong\prod_{Q(j,i)}M\]
    In particular, if $Q$ is an acyclic quiver, $i^\ast i_! M\cong M$ and $i^\ast i_\ast M\cong M$.
\end{rem}

Let $Q$ be a finite and acyclic quiver and, for any vertex $i\in Q_0$, define the $RQ$-module $U(i)$ such that $U(i)_k$ is equal to $R$ if $i=k$ and $0$ otherwise (so $U(i)_\alpha$ is the zero map for any $\alpha$). For a complex $X\in\cD(RQ)$, denote by $X(i)$ the tensor product $U(i)\boxtimes_{RQ}X$ and by $i_\times$ the functor $U(i)\boxtimes_{RQ}i_!:\cD(R)\to\cD(RQ)$. Notice that, by \cref{CompGen}, $U(i)\boxtimes_{RQ}\_$ preserves compact objects and so $i_\times$ does.

\begin{lem}\label{iprop}
    For any $X,Y\in\cD(RQ)$ and $M,N\in\cD(R)$ the following holds
    \begin{equation*}\begin{gathered}
        i_\times i^\ast X\cong X(i) \text{ and } j^\ast i_\times M\cong\begin{cases} M &\text{if }j=i\\\,0 &\text{otherwise}\end{cases} \\
        i^\ast(X\boxtimes_{RQ}Y)\cong i^\ast X\otimes_R i^\ast Y \text{ and } i_\times(M\otimes_R N)\cong i_\times M\boxtimes_{RQ} i_\times N
    \end{gathered}\end{equation*}
    In particular, we have that $i_\times(i^\ast X\otimes_R i^\ast Y)\cong X(i)\boxtimes_{RQ}Y(i)$.
\begin{proof}
    The first and the second equation easily follow from the definitions.
    
    As for the third one, note that the equality holds in the abelian setting, thus the left derived functors of the compositions $i^*(\_\boxtimes_{RQ}\_)$ and $i^*\_\otimes_Ri^*\_$ coincide. By the dual of \cite[Theorem 1.6.1]{Mil}, since $i^*$ is exact, the left derived functor of $i^*(\_\boxtimes_{RQ}\_)$ is the composition of the left derived functors of $i^*$ and $\boxtimes_{RQ}$. Moreover, since the adjunction $i^*\dashv i_*$ exists also between the homotopy categories (\cite[Proposition 1.1.3]{Mil}) and $i_*$ sends acyclic $R$-complexes to acyclic $RQ$-complexes, it follows that $i^*$ sends K-projective $RQ$-complexes to K-projective $R$-complexes and thus the dual of \cite[Theorem 1.6.1]{Mil} also applies to the left derived functor of $i^*\_\otimes_Ri^*\_$.
    
    The fourth equation follows from the fact that $i^*X\cong i^*Y$ implies $X(i)\cong Y(i)$ and the chain of quasi-isomorphisms $i^\ast i_!(M\otimes_R N)\cong i^\ast i_!M\otimes_R i^\ast i_!N\cong i^\ast(i_!M\boxtimes_R i_!N)$.
\end{proof}
\end{lem}

\section{The Balmer spectrum}\label{Sec2}
A \emph{thick tensor-ideal} of a tt-category $\cT$ is a thick subcategory $\cS$ such that $X\otimes S\in\cS$ for any $X\in\cT$ and $S\in\cS$. We will denote the smallest thick tensor-ideal containing a set of objects $\cX\subseteq\cT$ by $\langle\cX\rangle$, or by $\langle X\rangle$ when $\cX={X}$. A proper thick tensor-ideal $\cS$ is called \emph{prime} if whenever a product $X \otimes Y$ lies in $\cS$, then either $X\in\cS$ or $Y\in\cS$, or equivalently, if its complement $\cT\setminus\cS$ is $\otimes$-closed. The \emph{Balmer spectrum} of an essentially small tt-category $\cT$ is the set of the prime thick tensor-ideals of $\cT$, and it will be denoted by $\Spc(\cT)$. It has both a poset structure, given by inclusion, and a topological structure, given by the following basis of open subsets: $\leftcurly\{\cS\in\Spc(\cT)\mid X\in\cS\}\mid X\in\cT\rightcurly$. In this context, a \emph{specialization closed subset} of $\Spc(\cT)$ is equivalently a lower subset for the poset structure, or an arbitrary union of closed subsets for the topological structure. By \cite[Theorem 4.10, Remark 4.3, Remark 4.11]{Bal}, when $\Spc(\cT)$ is a noetherian topological space and $\cT$ is a rigid category, the specialization closed subsets of the Balmer spectrum fully classify the thick tensor-ideals of $\cT$.

Recall that the \emph{big support} of a complex $M\in\cD(R)$ is defined as
\[\Supp_R(M)=\{\p\in\Spec(R)\mid R_\p\otimesL M\neq0\}\]
and we denote by $K(\p)$ the Koszul complex at $\p$. The first part of next proposition follows similarly to \cite[Lemma 2.1.4.1]{LS}.

\begin{prop}\label{ThickFilt}
    For any finite and acyclic quiver $Q$ and $X\in\cD(RQ)$, it holds that
    \begin{enumerate}
        \item $\langle X\rangle=\thick_{RQ}\leftangle X(i)\mid i\in Q_0\rightangle$;
        \item If $X$ is compact, $\langle X\rangle=\thick_{RQ}\leftangle i_\times K(\p)\mid\p\in\Supp_R(i^\ast X),i\in Q_0\rightangle$.
    \end{enumerate}
    In particular, any thick tensor-ideal of $\cD^c(RQ)$ is generated by complexes of the form $i_\times K(\p)$ for some $i\in Q_0$ and $\p\in\Spec(R)$.
\begin{proof}
    (1) Since any finite acyclic quiver has a source and a sink, the $RQ$-module $U$ admits a filtration whose composition factors are isomorphic to the representations $U(i)$. Tensoring this filtration by a complex $X\in\cD(RQ)$ it follows that $X\in\thick_{RQ}\leftangle X(i)\mid i\in Q_0\rightangle$. Moreover, the latter is a tensor-ideal. Indeed, since thick subcategories of $\cD(R)$ are tensor-ideals, for any object $Z\in\cD(RQ)$, $i^\ast Z\otimes_R i^\ast X$ lies in the thick subcategory generated by $i^\ast X$. Thus, since $Z\boxtimes_{RQ}X(i)=Z(i)\boxtimes_{RQ}X(i)$, by \cref{iprop}, $Z\boxtimes_{RQ}X(i)=i_\times(i^\ast Z\otimes_R i^\ast X)$ which lies in the thick subcategory generated by $i_\times i^\ast X$, i.e.~$\thick_{RQ}\leftangle X(i)\rightangle$. On the other hand, by the closure under tensor product, we have $X(i)=U(i)\boxtimes_{RQ}X\in\leftangle X\rightangle$ for any $i\in Q_0$.

    \noindent(2) If $X\in\cD^c(RQ)$, by \cite[Theorem 1.5]{NeeTC}, we have that the thick subcategory $\thick_R\langle i^\ast X\rangle$ is equal to $\thick_R\leftangle K(\p)\mid\p\in\Supp_R(i^\ast X)\rightangle$. Then, applying the functor $i_\times$, by \cref{iprop}, we can conclude that
    \[\thick_{RQ}\langle X(i)\rangle=\thick_{RQ}\leftangle i_\times K(\p)\mid\p\in\Supp_R(i^\ast X)\rightangle\]
    Thus, combining this with point (1) gives the statement.

    \noindent As for the last part, it is sufficient to note that, given a set of objects $\cX\subseteq\cD^c(RQ)$, the thick tensor-ideal $\langle\cX\rangle$ is equal to $\thick_{RQ}\leftangle i_\times K(\p)\mid\p\in\Supp_R(i^\ast X),X\in\cX,i\in Q_0\rightangle$.
\end{proof}
\end{prop}

For any prime ideal $\p\in\Spec(R)$ and $i\in Q_0$, define the functor $\xi_{\p,i}:=R_\p\otimes_R i^\ast(\_)$ such that
\[\xi_{\p,i}:\cD^c(RQ)\to\cD^c(R_\p)\]
\[X\mapsto i^\ast X_\p\]
Denote by $\cS_{\p,i}:=\ker\xi_{\p,i}=\{X\in\cD^c(RQ)\mid i^*X_\p=0\}$. Note that, since the functor $\xi_{\p,i}$ is exact and commutes with $\boxtimes_{RQ}$, each subcategory $\cS_{\p,i}$ is a thick tensor-ideal. In particular, by \cref{ThickFilt} (2):
\[\cS_{\p,i}=\thick_{RQ}\leftangle\{i_\times K(\q)\mid\q\nsubseteq\p\}\cup\{j_\times K(\q)\mid j\neq i,\q\in\Spec(R)\}\rightangle\]

In \cite{LS}, the authors prove that for any field $\bK$, the Balmer spectrum of $(\cD(\bK Q), \boxtimes_{\bK Q})$ is a discrete space consisting of as many points as the vertices of the quiver $Q$. In the following theorem, we generalize \cite[Theorem 2.1.5.1]{LS} to any commutative noetherian ring, showing that in our context the Balmer spectrum is a disjoint union of as many copies of the prime spectrum $\Spec(R)$ as the number of vertices of $Q$. Moreover, the following theorem also provides a direct proof, when $Q=A_1$, of the homeomorphism between the Balmer spectrum of $\cD^c(R)$ and the prime spectrum of the ring $\Spec(R)$, which originally follows by \cite[Corollary 5.6]{Bal}.

\begin{thm}\label{Spectrum}
    For any commutative noetherian ring $R$ and finite acyclic quiver $Q$, we have that
    \begin{enumerate}
        \item $\Spc(\cD^c(RQ))=\{\cS_{\p,i}\mid\p\in\Spec(R),i\in Q_0\}$;
        \item There is an order-reversing homeomorphism
        \[f:\Spc(\cD^c(RQ))\longrightarrow\Spec(R)\times Q_0\]
        i.e.~$\Spc(\cD^c(RQ))\overset{f}{\cong}_\mathrm{\,Top}\Spec(R)\times Q_0$ and $\Spc(\cD^c(RQ))\overset{f}{\cong}_{\,\Pos}\Spec(R)^\mathrm{op}\times Q_0$.
    \end{enumerate}
\begin{proof}
    (1) Let us start by checking that each thick tensor-ideal $\cS_{\p,i}$ is prime by showing that their complements are $\boxtimes$-closed. Given two complexes $X,Y\in\cD^c(RQ)\setminus\cS_{\p,i}$, since $\cS_{\p,i}=\ker\xi_{\p,i}$, the complexes $i^\ast X_\p$ and $i^\ast Y_\p$ are non-zero, bounded and vertexwise finitely generated $R_\p$-free. In particular, by the so-called K{\"u}nneth formula (see \cite[2.5.18 (c)]{CFH}), the top cohomology of $i^\ast X_\p\otimes_R i^\ast Y_\p$ is the tensor product of the top cohomologies of $i^\ast X_\p$ and $i^\ast Y_\p$ and so, by the Nakayama's lemma (see \cite[Exercise 2.3]{AM}) it is non-zero. Thus, by \cref{iprop}, $i^\ast (X\boxtimes_{RQ}Y)_\p\neq0$, i.e.~$X\boxtimes_{RQ}Y$ does not lie in $\cS_{\p,i}$. Now let us prove that all thick tensor-ideal are of this form. Let $\cS$ be a prime thick tensor-ideal of $\cD^c(RQ)$, by \cref{ThickFilt} (2) it is generated by a proper subset of $\cK=\leftcurly i_\times K(\p)\mid\p\in\Spec(R),i\in Q_0\rightcurly$. Let $\cS^\complement=\cD^c(RQ)\setminus\cS$ be the complement of $\cS$ and consider the set $\cK\cap\cS^\complement$. Since $\cS$ is proper, there exists a vertex $\bar{i}\in Q_0$ such that $U\left(\bar{i}\right)\in\cS^\complement$ (otherwise $U\in\cS$ and $\cS=\cD^c(RQ)$). Since $\cS^\complement$ is $\boxtimes$-closed and $0\notin\cS^\complement$, we have that $\cK\cap\cS^\complement\subseteq\leftcurly \bar{i}_\times K(\p)\mid\p\in\Spec(R)\rightcurly$ (indeed for any $j\neq\bar{i}$ and $\q\in\Spec(R)$, we have $j_\times K(\q)\boxtimes_{RQ}\bar{i}_\times K(\p)=0$). Moreover, since $R$ is noetherian, the set $\{\p\in\Spec(R)\mid \bar{i}_\times K(\p)\in\cS^\complement\}$ has a maximal element, say $\bar{\p}$, and for any $\q\subseteq\bar{\p}$ the complex $\bar{i}_\times K(\q)\in\cS^\complement$ (otherwise, by \cite[Lemma 1.2]{NeeTC}, $K\left(\bar{\p}\right)\in\thick_R\langle K(\q)\rangle$ and so $\bar{i}_\times K\left(\bar{\p}\right)\in\cS$). We end by proving that the maximal element $\bar{\p}$ is unique and so conclude that $\cS=\cS_{\bar{\p},\bar{i}}$. Suppose that there is another maximal element $\bar{\q}$, then:
    \begin{itemize}
        \item[(i)\,] If $\bar{\p}+\bar{\q}=R$, we have $K\left(\bar{\p}\right)\otimes_R K\left(\bar{\q}\right)=0$. Indeed, by \cite[14.3.1]{CFH}, the tensor product is equal to $K\left(\bar{\p}+\bar{\q}\right)$ which is an acyclic complex by \cite[11.4.17]{CFH}. In particular, the tensor product $\bar{i}_\times K\left(\bar{\p}\right)\boxtimes_{RQ} \bar{i}_\times K\left(\bar{\q}\right)=0$ lies in $\cS$;
        \item[(ii)] If $\bar{\p}+\bar{\q}\neq R$, we have $K\left(\bar{\p}\right)\otimes_R K\left(\bar{\q}\right)=K\left(\bar{\p}+\bar{\q}\right)$ by \cite[14.3.1]{CFH}. Let us prove that, even in this case, the tensor product $\bar{i}_\times K\left(\bar{\p}\right) \boxtimes_{RQ}\bar{i}_\times K\left(\bar{\q}\right)=\bar{i}_\times K\left(\bar{\p}+\bar{\q}\right)$ lies in $\cS$. Indeed, there exists finitely many prime ideals $\m_1,\ldots,\m_N$ such that $\bar{\p}+\bar{\q}\subseteq\m_k$ for any $k=1,\ldots,N$ and $V\left(\bar{\p}+\bar{\q}\right)=V(\m_1)\cup\ldots\cup V(\m_N)$ (see \cite[Lemma 10.63.5-6]{STA}). So, by maximality of $\bar{\p}$, it holds that $\bar{i}_\times K(\m_k)\in\cS$ for any $k=1,\ldots,N$ and, by \cite[Lemma 1.2]{NeeTC}, the complex $K\left(\bar{\p}+\bar{\q}\right)$ lies in $\thick_R\langle K(\m_1)\oplus\ldots\oplus K(\m_N)\rangle$. Thus, it follows that $\bar{i}_\times K\left(\bar{\p}+\bar{\q}\right) \in\cS$.
    \end{itemize}
    Since $\cS^\complement$ is $\boxtimes$-closed, both cases lead to a contradiction. Thus, $\cK\cap\cS^\complement=\leftcurly\bar{i}_\times K(\q)\mid\q\subseteq\bar{\p} \rightcurly$ and so $\cS=\cS_{\bar{\p},\bar{i}}$.

    \noindent(2) Recall that the topology on $\Spc(\cD^c(RQ))$ is given by the basis of open subsets $\{\{\cS_{\p,i}\mid X\in\cS_{\p,i}\}\mid X\in\cD^c(RQ)\}$, while the topology on $\Spec(R)\times Q_0$, where $Q_0$ is taken with discrete topology, is given by the product topology, thus by the basis of open subsets $\{\{\p\mid r\notin\p\}\times\{i\}\mid r\in R,i\in Q_0\}$. We want to prove that the bijection given by the assignment $f:\cS_{\p,i}\longmapsto(\p,i)$ is a homeomorphism. Let us start proving the continuity. Given an open set of the basis in the codomain $\cO=\{\p\mid r\notin\p\}\times\{i\}$, for an element of the ring $r\in R$, by \cite[15.1.10]{CFH}, we have that $\Supp_R(K((r)))=V((r))$. Thus, it holds that
    \[\xi_{\q,j}\left(i_\times K((r))\oplus\left(\textstyle{\bigoplus}_{\ell\neq i}\,\ell_\times R\right)\right)=0 \text{ if and only if } j=i \text{ and } r\notin\q\]
    this implies $f^{-1}(\cO)=\leftcurly\cS_{\p,i}\mid i_\times K((r))\oplus\left(\bigoplus_{\ell\neq i}\,\ell_\times R\right)\in\cS_{\p,i}\rightcurly$ which is an open set of the basis.
    
    Moreover, $f$ is also an open map. Indeed, given an open set of the basis in the domain $\cO'=\{\cS_{\p,i}\mid X\in\cS_{\p,i}\}$, for some complex $X\in\cD^c(RQ)$, we have that
    \[f(\cO')=\{(\p,i)\mid\xi_{\p,i}(X)=0\}=\bigcup_{i\in Q_0}\Supp_R(i^*X)^\complement\times\{i\}\]
    Since the functor $i^*$ preserves compacts and the support of a compact object of $\cD(R)$ is closed in $\Spec(R)$ (see \cite[14.1.6]{CFH}), it follows that $f(\cO')$ is an open set.
    As for the order-reversing property, note that for any $\cS_{\q,j}\subseteq\cS_{\p,i}$, we have that $\cK\cap\cS_{\p,i}^\complement\subseteq \cK\cap\cS_{\q,j}^\complement$. In particular, by the previous point, it implies that $i=j$ and $\p\subseteq\q$, i.e.~$(\p,i)\leq(\q,j)$.
\end{proof}
\end{thm}

\begin{exmp}
    Let $A_2=\overset{1}{\bullet}\longrightarrow\overset{2}{\bullet}\ $ and consider the algebra $\bK[[x]]A_2$, then we have the posets
    \[\Spc(\cD^c(\bK[[x]]A_2))=\begin{tikzcd} {\cS_{(0),1}} \\ \\ {\cS_{(x),1}} \arrow[uu] \end{tikzcd} \ \textstyle{\bigcup}\ \begin{tikzcd} {\cS_{(0),2}} \\ \\ {\cS_{(x),2}} \arrow[uu] \end{tikzcd} \qquad\qquad \Spec(\bK[[x]])\times\{1,2\}=\begin{tikzcd} {((x),1)} \\ \\
    {((0),1)} \arrow[uu] \end{tikzcd} \ \textstyle{\bigcup}\ \begin{tikzcd} {((x),2)} \\ \\ {((0),2)} \arrow[uu] \end{tikzcd}\]
    In this case, we have that $K((x))\cong\bK$, thus the open set of the basis $\{\p\mid x\notin\p\}\times\{1\}=\{((0),1)\}$ correspond to the open set of he basis $\{\cS_{\p,i}\mid 1_\times\bK\oplus2_\times\bK[[x]]\in\cS_{\p,i}\}=\{\cS_{(0),1}\}$.
\end{exmp}

\subsection{Thick tensor-ideals}
We saw in \cref{notrigid} that the tt-category $(\cD^c(RQ),\boxtimes_{RQ},\mathbf{U})$ is not rigid. In any case, rigidity is not a necessary condition to get a full classification of all the thick tensor-ideals. For this it suffices that all the thick tensor-ideals are radical.

\begin{prop}\label{Rigid}
    For any $X\in\cD^c(RQ)$, it holds that $X\in\leftangle X\boxtimes_{RQ}X\rightangle$. In particular, any thick tensor-ideal is radical.
\begin{proof}
    First note that the results holds in $\cD^c(R)$, indeed there the tensor products $\otimes$ and $\boxtimes$ coincide and any object is $\otimes$-rigid. So, for any $X\in\cD^c(RQ)$, we have that $\thick_R\leftangle i^\ast X\rightangle=\thick_R\leftangle i^\ast X\otimes_R i^\ast X\rightangle$ and so $\thick_{RQ}\leftangle i_\times i^\ast X\rightangle=\thick_{RQ}\leftangle i_\times(i^\ast X\otimes_R i^\ast X)\rightangle$. Thus, by \cref{iprop}, it follows that $\thick_{RQ}\leftangle X(i)\rightangle=\thick_{RQ}\leftangle X(i)\boxtimes_{RQ}X(i)\rightangle$. Combining this with \cref{ThickFilt} (1) we get the statement and, by \cite[Proposition 4.4]{Bal}, it follows that any thick tensor-ideal in $\cD^c(RQ)$ is radical.
\end{proof}
\end{prop}

It follows the classification of the thick tensor-ideals of $(\cD^c(RQ),\boxtimes_{RQ},\mathbf{U})$. Here we get the result from the general theory of tt-categories in \cite{Bal}.

\begin{cor}\label{Thickbyspc}
    For any commutative noetherian ring $R$ and finite acyclic quiver $Q$, it holds that
    \[\Thick_\boxtimes(\cD^c(RQ))\longleftrightarrow\cV(\Spc(\cD^c(RQ)))\]
    where $\cV(\Spc(\cD^c(RQ)))$ is the set of specialization closed subsets of $\Spc(\cD^c(RQ))$.
\begin{proof}
    By \cref{Spectrum}, the Balmer spectrum $\Spc(\cD^c(RQ))$ is a finite union of noetherian spaces, indeed $\Spec(R)\times Q_0\cong\bigsqcup_{Q_0}\Spec(R)$, and so it is noetherian too. Thus, by \cite[Theorem 4.10, Remark 4.11]{Bal} and \cref{Rigid}, the thick tensor-ideals of $\cD^c(RQ)$ are in bijection with the specialization closed subsets of $\Spc(\cD^c(RQ))$.
\end{proof}
\end{cor}

We can also give two alternative classifications of the thick tensor-ideals of $\cD^c(RQ)$, following from \cref{Thickbyspc}. Here, we present a proof using the standard adjunction isomorphism in the category of posets:
\[\Hom_{\Pos}(X \times A,B)\cong\Hom_{\Pos}(A,\Hom_{\Pos}(X,B)) \text{ for any three posets } X,A \text{ and } B\]
given by the assignment $\phi : f \mapsto \phi_f$, where $\phi_f(a)(x) = f((x,a))$, with inverse $\psi : g \mapsto \psi_g$, where $\psi_g((x,a)) = g(a)(x)$. A direct proof also follows from the arguments used in \cref{Aislebyaisle}.

\begin{cor}\label{Thickbythick}
    For any commutative noetherian ring $R$ and finite acyclic quiver $Q$, it holds that
    \begin{enumerate}
        \item $\Thick_\boxtimes(\cD^c(RQ))\cong\Thick(\cD^c(R))\times Q_0$;
        \item $\Thick_\boxtimes(\cD^c(RQ))\cong\Hom_{\Pos}(\Spec(R),\Thick_\boxtimes(\cD(\bK Q)))$
    \end{enumerate}
\begin{proof}
    (1) Let $\{0,1\}$ be the poset with two elements and the obvious order, then the specialization closed subsets of $\Spc(\cD^c(RQ))$ can be identified with $\Hom_{\Pos}(\Spc(\cD^c(RQ)),\{0,1\})$ by associating to each subset its characteristic function. From the standard adjunction isomorphism in the category of posets, we get that
    \[\Hom_{\Pos}(\Spec(R)\times Q_0,\{0,1\})\cong\Hom_{\Pos}(Q_0,\Hom_{\Pos}(\Spec(R),\{0,1\}))\]
    where, by \cite[Theorem 1.5]{NeeTC}, $\Hom_{\Pos}(\Spec(R),\{0,1\}))\cong\Thick(\cD^c(R))$ and thus, since $Q_0$ is a discrete poset, we get the statement.

    \noindent (2) Alternatively, applying the other possible adjunction isomorphism we get
    \[\Hom_{\Pos}(\Spec(R)\times Q_0,\{0,1\})\cong\Hom_{\Pos}(\Spec(R),\Hom_{\Pos}(Q_0,\{0,1\}))\]
    where, since $Q_0$ is a discrete poset, $\Hom_{\Pos}(Q_0,\{0,1\})$ is $\cP(Q_0)$, the power set of $Q_0$. In virtue of the bijection between Serre subcategories and sets of simple modules (\cite[Section 8]{Kan}), $\cP(Q_0)$ is isomorphic to the lattice $\Serre(\bK Q)$ and thus by \cref{OverField} we get the statement.
\end{proof}
\end{cor}

\begin{rem}
    Comparing \cref{Thickbythick} (2) with the classification in \cite{AS}, where the authors prove that $\Thick(\cD^c(RQ))\cong\Hom_{\Pos}(\Spec(R),\mathbf{Nc}(Q))$, one can note that the gap between the thick tensor-ideals and thick subcategories in $\cD^c(RQ)$ is reflected by the gap between the lattices $\cP(Q_0)$ and $\mathbf{Nc}(Q)$, the lattice of noncrossing partitions of $Q$.
\end{rem}

\subsection{Stratification}\label{StratSec}
As seen in \cref{Thickbyspc}, although the tt-category $(\cD^c(RQ), \boxtimes_{RQ}, \mathbf{U})$ is almost never rigid, we can still recover the classification of all thick tensor-ideals from its Balmer spectrum. However, without rigidity we loose many aspects of the Balmer-Favi theory in \cite{BF}, which describes how the Balmer spectrum of the compact objects, $\Spc(\cT^c)$, controls the ambient compactly generated tt-category $\cT$.

In a rigidly-compactly generated tt-category $(\cT,\otimes,\mathbf{1})$, any smashing tensor-ideal $\cS$ gives rise to a tensor-idempotent triangle, i.e.~a distinguished triangle
\[\begin{tikzcd} e(\cS) \arrow[r] & \mathbf{1} \arrow[r] & f(\cS) \arrow[r,"+"] & {} \end{tikzcd}\]
such that the following (equivalent by \cite[Proposition 3.1]{BF}) identities hold
\[e(\cS) \otimes e(\cS) = e(\cS) \qquad e(\cS) \otimes f(\cS) = 0 \qquad f(\cS) \otimes f(\cS) = f(\cS)\]
Moreover, by \cite[Theorem 3.5]{BF}, in this context, there is a bijection between smashing tensor-ideals and tensor-idempotent triangles. The key point is that, by \cite[Theorem 2.13]{BF} and \cite[Definition 3.3.2]{HPS}, a tensor-ideal $\cS$ is smashing if and only if its right orthogonal $\cS^\perp$ is also a tensor-ideal -- a property that crucially depends on rigidity.

The next example show that, for a general compactly generated tt-category, we loose the aforementioned bijection and that there can be more smashing tensor-ideals than tensor-idempotent triangles.

\begin{exmp}\label{MoreEx}
    Let $A_2=\overset{1}{\bullet}\longrightarrow\overset{2}{\bullet}\ $ and consider the algebra $\bK A_2$. Then the lattice of smashing tensor-ideals is
    
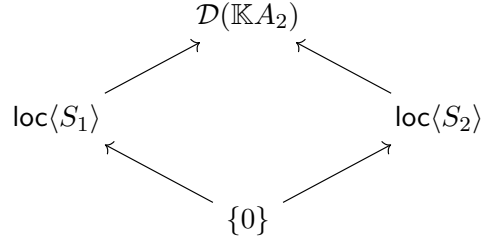
\[\begin{tikzcd}
         & \cD(\bK A_2) & \\
        \loc\langle S_1\rangle \arrow[ru] & & \loc\langle S_2\rangle \arrow[lu] \\
         & \{0\} \arrow[lu] \arrow[ru] &
    \end{tikzcd}\]
    \noindent and the truncation triangles of $\mathbf{U}=\bK\xrightarrow{1}\bK$ with respect to the tensor-ideals $\loc\langle S_1\rangle$ and $\{0\}$ are both $0\to\mathbf{U}\to\mathbf{U}\xrightarrow{+}\ $. While the ones with respect to $\loc\langle S_2\rangle$ and $\cD(\bK A_2)$ are $S_2\to\mathbf{U}\to S_1\xrightarrow{+}\ $ and $\mathbf{U}\to\mathbf{U}\to0\xrightarrow{+}\ $, respectively.
\end{exmp}

This implies that, a priori, \emph{tensor-idempotent residue objects} objects may not exist for all points of the spectrum (\cite[Definition 7.9]{BF}), and consequently, the \emph{Balmer-Favi support} for non-compact objects may not be defined (\cite[Definition 7.16]{BF}). Indeed, in a rigidly-compactly generated tt-category $(\cT,\otimes,\mathbf{1})$, for any prime thick tensor-ideal $\cP$ in $\Spc(\cT^c)$, we can find two smashing tensor-ideals, $\cV_\cP$ and $\cZ_\cP$, and define the tensor-idempotent residue object at $\cP$ as
\[g(\cP):=e(\cV_\cP)\otimes f(\cZ_\cP)\]
The Balmer-Favi support of an object $X \in \cT$ is then given by
\[\Supp_\cT(X)=\{\cP\in\Spc(\cT^c)\mid X\otimes g(\cP)\neq 0\}\]
As shown in \cite[Proposition 7.18]{BF}, the tensor-idempotent residue object $g(\cP)$ uniquely determines the prime $\cP$, in the sense that $\Supp_\cT(g(\cP))=\{\cP\}$.

The major implication of the absence of these notions is the lack of a stratification theory for the classification of localizing tensor-ideals (see \cite[Section 4]{BHS}). Following \cite[Theorem A]{BHS}, a rigidly-compactly generated tt-category $\cT$ is said to be \emph{stratified} (by the Balmer-Favi support) if, denoting by $\langle X\rangle^\amalg$ the smallest localizing tensor-ideal generated by an object $X\in\cT$, the equivalent conditions hold:
\begin{itemize}
    \item[$\boldsymbol{\cdot}$](\emph{Local-to-global and minimality}) For any object $X\in\cT$, $X\in\leftangle g(\cP)\otimes X\mid\cP\in\Spc(\cT^c)\rightangle^\amalg$ and $\leftangle g(\cP)\rightangle^\amalg$ is a minimal localizing tensor-ideal for any $\cP\in\Spc(\cT^c)$;
    \item[$\boldsymbol{\cdot}$](\emph{Stratification}) Taking the Balmer-Favi support and its preimage defines a bijection
    \[\Loc_\otimes(\cT)\longleftrightarrow\leftcurly\text{Subsets of }\Spc(\cT^c)\rightcurly\]
\end{itemize}

Nevertheless, even in the absence of a general theory, it is still possible -- based on the proof of \cref{ThickFilt} and the local-to-global principle in $\cD(R)$ (\cite[Section 3.4]{Stour}) -- to construct \emph{ad hoc} residue objects $g(\p,i)$ for all thick prime tensor-ideals $\cS_{\p,i}$ and thereby obtain a stratification result (\cref{Strat}).

\begin{defn}\label{BigSupp}
    For a prime thick tensor-ideal $\cS_{\p,i}$ in $\Spc(\cD^c(RQ))$, we define the tensor-idempotent residue object at $\cS_{\p,i}$ to be the complex
    \[g(\p,i):=i_\times K_\infty(\p)_\p\]
    where $K_\infty(\p)$ is the infinite Koszul complex at $\p$. Consequently, we define the support of an object $X\in\cD(RQ)$ as
    \[\Supp_{RQ}(X)=\{\cS_{\p,i}\in\Spc(\cD^c(RQ))\mid g(\p,i)\boxtimes_{RQ}X\neq0\}\] 
\end{defn}

\begin{rem}
    It is not clear how to recover this notion of support for a general compactly generated tt-category. Indeed, as we can see in \cref{MoreEx}, the residue objects do not arise, in any reasonable way, from truncation triangles of the unit $\mathbf{U}$ with respect to some smashing tensor-ideals.
\end{rem}

The conditions of \cite[Theorem A]{BHS} reported above are satisfied by the tt-category $\cD(RQ)$.
\pagebreak

\begin{thm}\label{Strat}
    For any finite and acyclic quiver $Q$, the following statements hold:
    \begin{enumerate}
    \item For any $X\in\cD(RQ)$, $\langle X\rangle^\amalg=\leftangle g(\p,i)\boxtimes_{RQ}X\mid\cS_{\p,i}\in\Spc(\cD^c(RQ))\rightangle^\amalg$ and $\leftangle g(\p,i)\rightangle^\amalg$ is a minimal localizing tensor-ideal for any $\cS_{\p,i}\in\Spc(\cD^c(RQ))$;
    \item Taking the support and its preimage defines a bijection
    \[\Loc_\boxtimes(\cD(RQ))\longleftrightarrow\leftcurly\text{Subsets of }\Spc(\cD^c(RQ))\rightcurly\]
\end{enumerate}
\begin{proof}
    (1) For any $X\in\cD(RQ)$, by \cref{ThickFilt} (1), we have that
    \[\langle X\rangle^\amalg=\loc_{RQ}\langle X(i)\mid i\in Q_0\rangle\tag{$\ast$}\]
    Moreover, by the local-to-global principle in $\cD(R)$ (\cite[Section 3.4]{Stour}), for any $i\in Q_0$, it holds that
    \[\loc_R\leftangle i^\ast X\rightangle=\loc_R\leftangle K_\infty(\p)_\p\otimes_R i^\ast X\mid \p\in\Spec(R)\rightangle\]
    So, applying the functor $i_\times$, by \cref{iprop}, we can conclude that
    \[\loc_{RQ}\langle X(i)\mid i\in Q_0\rangle=\loc_{RQ}\langle i_\times K_\infty(\p)_\p\boxtimes_{RQ}X(i)\mid\cS_{\p,i}\in\Spc(\cD^c(RQ))\rangle\]
    Noting that $i_\times K_\infty(\p)_\p\boxtimes_{RQ}X(i)=(g(\p,i)\boxtimes_{RQ}X)(i)$, by $(\ast)$ the above implies that
    \[\langle X\rangle^\amalg=\langle g(\p,i)\boxtimes_{RQ}X\mid\cS_{\p,i}\in\Spc(\cD^c(RQ))\rangle^\amalg\]
    Moreover, by \cite[3.24]{Stour}, for any $\p\in\Spec(R)$ the subcategory $\loc_R\leftangle K_\infty(\p)_\p\rightangle$ is minimal in $\cD(R)$. Then, for any localizing tensor-ideal $\cL$ of $\cD(RQ)$, the inclusion $\cL\in\leftangle g(\p,i)\rightangle^\amalg$ implies both $\cL=i_\times i^\ast(\cL)$ and $i^\ast(\cL)\subseteq\loc_R\leftangle K_\infty(\p)_\p\rightangle$, respectively by $(\ast)$ and \cref{iprop}. Thus, we can deduce the minimality of the localizing tensor-ideal $\leftangle g(\p,i)\rightangle^\amalg$.

    \noindent (2) For any object $X\in\cD(RQ)$, by minimality, we have that $\leftangle g(\p,i)\boxtimes_{RQ}X\rightangle^\amalg=\leftangle g(\p,i)\rightangle^\amalg$ if $S_{\p,i}\in\Supp_{RQ}(X)$ or zero otherwise. Defining the support of a subcategory $\cL\subseteq\cD(RQ)$ as the union of the supports of its objects, from the previous point, we have that any localizing tensor-ideal $\cL$ is equal to $\leftangle g(\p,i)\mid\cS_{\p,i}\in\Supp_{RQ}(\cL)\rightangle$. Noting that the $\Supp_{RQ}(g(\p,i))=\{\cS_{\p,i}\}$, we have the bijection.
\end{proof}
\end{thm}

\section{Tensor-t-structures}\label{Sec3}
We recall the definition of tensor-t-structure from \cite{DS}. Let $(\cT,\otimes,\mathbf{1})$ be a compactly generated tt-category and suppose it is given with a suspended subcategory $\cT^{\leq 0}$ such that $\cT^{\leq 0}\otimes\cT^{\leq 0}\subseteq\cT^{\leq 0}$ and $\mathbf{1}\in\cT^{\leq 0}$. A suspended subcategory $\cS$ of $\cT$ is called \emph{tensor-suspended} (\emph{with respect to $\cT^{\leq 0}$}) if $\cT^{\leq 0}\otimes\cS\subseteq\cS$. A t-structure $(\cU,\cV)$ of $\cT$ is called a \emph{tensor-t-structure} (with respect to $\cT^{\leq 0}$) if the aisle $\cU$ is a tensor-suspended subcategory -- or equivalently if $[\cT^{\leq 0},\cV]\subseteq\cV$ (see \cite[Proposition 3.9]{DS}). In this case $\cU$ and $\cV$ are called \emph{tensor-aisle} and \emph{tensor-coaisle}, respectively. Given a set of objects $\cX\subseteq\cT$, we denote the smallest cocomplete tensor-suspended subcategory containing $\cX$ by $\langle\cX\rangle^\leq$ and we call a tensor-aisle $\cU$ \emph{compactly generated} if $\cU=\langle\cX\rangle^\leq$ for a set of compact objects $\cX\subseteq\cT^c$. Recall that, by \cite[Proposition 3.11 (i)]{DS}, whenever $\cT^{\leq 0}$ is compactly generated, $\langle\cS\rangle^\leq$ is a tensor-aisle for any set of compact objects $\cS\subseteq\cT^c$ and, moreover, note that this result holds for general compactly generated tt-categories, since the rigidity assumption do not play any role.

In the following we will consider tensor-t-structures of $\cD(RQ)$ with respect to the aisle of the standard t-structure $(\cD^{\leq 0}(RQ),\cD^{\geq 1}(RQ))$. Observe that, indeed, $\mathbf{U}\in\cD^{\leq 0}(RQ)$ and $\cD^{\leq 0}(RQ)\boxtimes_{RQ}\cD^{\leq 0}(RQ)\subseteq\cD^{\leq 0}(RQ)$. We denote by $\susp_{RQ}^\amalg\leftangle\cX\rightangle$ the smallest cocomplete suspended subcategory of $\cD(RQ)$ containing a set of object $\cX$ and recall that, by \cite[Theorem 3.4]{AJSo}, it is the aisle of a t-structure. Analogously to \cref{ThickFilt}, we have the following.

\begin{prop}\label{AisleFilt}
    For any finite and acyclic quiver $Q$ and $X\in\cD(RQ)$, it holds that:
    \begin{enumerate}
        \item $\langle X\rangle^\leq=\susp_{RQ}^\amalg\leftangle X(i)\mid i\in Q_0\rightangle$;
        \item If $X$ is compact, $\langle X\rangle^\leq=\susp_{RQ}^\amalg\leftangle i_\times K(\p)[-n]\mid\p\in\Supp_R(H^n(i^\ast X)),n\in\bZ,i\in Q_0\rightangle$.
    \end{enumerate}
    In particular, any cocomplete tensor-suspended subcategory of $\cD^c(RQ)$ is generated by complexes of the form $i_\times K(\p)[-n]$ for some $i\in Q_0$, $\p\in\Spec(R)$ and $n\in\bZ$.
\begin{proof}
    (1) Noting that all the $RQ$-modules $U(i)$ belong to the standard aisle $\cD^{\leq 0}(RQ)$, the proof follows as the one of \cref{ThickFilt} (1), since there we do not use closure under negative shifts.

    \noindent(2) If $X\in\cD^c(RQ)$, by the main result in \cite{AJS}, we have that
    \[\susp_R^\amalg\langle i^\ast X\rangle=\susp_R^\amalg\leftangle K(\p)[-n]\mid\p\in\Supp_R(H^n(i^\ast X)),n\in\bZ\rightangle\]
    So, applying the functor $i_\times$, by \cref{iprop}, we can conclude that
    \[\susp_R^\amalg\langle X(i)\rangle=\susp_R^\amalg\leftangle i_\times K(\p)[-n]\mid\p\in\Supp_R(H^n(i^\ast X)),n\in\bZ\rightangle\]
    and combining this with point (1) gives the statement.

    \noindent As for the last part, it is sufficient to note that, given a set of objects $\cX\subseteq\cD^c(RQ)$, the cocomplete tensor-suspended subcategory $\langle\cX\rangle^\leq$ is equal to
    \[\susp_{RQ}^\amalg\leftangle i_\times K(\p)[-n]\mid\p\in\Supp_R(H^n(i^\ast X)),X\in\cX,n\in\bZ,i\in Q_0\rightangle\]
\end{proof}
\end{prop}

For a tt-category $(\cT,\otimes,\mathbf{1})$, we denote by $\Susp_\otimes(\cT^c)$ the lattice of tensor-suspended subcategories of $\cT^c$ ordered by inclusion and by $\mathrm{Aisle_{cg}}_\otimes(\cT)$ the lattice of compactly generated tensor-aisle of $\cT$. Recall that, by \cite[Theorem 15 (i)]{DSw}, these two lattices are isomorphic and note that this result holds for compactly generated tt-categories without the rigidity assumption.

Now we present the analogous of \cref{Thickbythick} for compactly generated tensor-aisles of $\cD(RQ)$. We start by proving a version of \cite[Theorem A.3]{Sab} adapted to tensor-aisles.

\begin{lem}\label{OverField}
    For any field $\bK$ and finite acyclic quiver $Q$, we have that
    \[\Aisle_\boxtimes(\cD(\bK Q))\cong\Filt(\Serre(\bK Q))\]
    In particular, the thick tensor-ideals of $\cD^c(\bK Q)$ are in bijection with $\Serre(\bK Q)$.
\begin{proof}
    Since any aisle in $\cD(\bK Q)$ is compactly generated, by \cite[Theorem 15 (i)]{DSw}, it suffices to prove that the following bijection holds
    \begin{equation*}\begin{gathered}
        \Susp_\boxtimes(\cD^c(\bK Q))\longleftrightarrow\Filt(\Serre(\bK Q)) \\
        \qquad\qquad\cS\longmapsto\left(\ldots\supseteq H^n(\cS)\supseteq H^{n+1}(\cS)\supseteq\ldots\right) \\
        \leftcurly X\in\cD^c(\bK Q)\mid H^n(X)\in\cZ_n\rightcurly\longmapsfrom\,\left(\ldots\supseteq\cZ_n\supseteq\cZ_{n+1}\supseteq\ldots\right)
    \end{gathered}\end{equation*}
    Given a tensor-suspended subcategory $\cS\subseteq\cD^c(\bK Q)$, by \cite[Lemma A.1]{Sab}, the class $H^n(\cS)$ is closed under images, cokernels and extensions for any $n\in\bZ$. So, it is sufficient to prove that it is closed also under subobjects to prove that it is a Serre subcategory of $\mods(\bK Q)$. Let $M\in H^n(\cS)$, then $M[-n]$ lies in $\cS$ and, by closure under tensor product, all the modules of the form $U(i)\boxtimes_{\bK Q}M[-n]$ belongs to $\cS$. In particular, all the composition factors of $M$ lies in $H^n(\cS)$ and so does any of its submodules. Moreover, for a filtration of Serre subcategories $\ldots\supseteq\cZ_n\supseteq\cZ_{n+1}\supseteq\ldots$, the class $\cS_\cZ=\{X\in\cD^c(\bK Q)\mid H^n(X)\in\cZ_n\}$ is clearly closed under summands and positive shifts, so it remains to show that it is closed under extensions. Let $X,Y\in\cS_\cZ$ and consider a distinguished triangle
    \[\begin{tikzcd} X \arrow[r] & Z \arrow[r] & Y \arrow[r,"+"] & {} \end{tikzcd}\]
    It induces the exact sequence in cohomology
    \[\begin{tikzcd}
    H^{n-1}(Y) \arrow[r,"f"] & H^n(X) \arrow[r] & H^n(Z) \arrow[r] & H^n(Y) \arrow[r,"g"] & H^{n+1}(X)
    \end{tikzcd}\]
    from which we obtain the short exact sequence
    \[\begin{tikzcd} 0 \arrow[r] & \coker f \arrow[r]  & H^n(Z) \arrow[r] & \ker g \arrow[r] & 0 \end{tikzcd}\]
    Since both $\coker f$ and $\ker f$ lies in $\cZ_n$ so does $H^n(Z)$. Since $\bK Q$ is hereditary, it follows that $\cS_\cZ=\susp_{\bK Q}^\amalg\langle\cZ_n[-n]\mid n\in\bZ\rangle$ and, since any Serre subcategory $\cZ_n$ is uniquely determined by the simple modules it contains (\cite[Section 8]{Kan}), it is actually generated by shifts of simple $\bK Q$-modules. In particular, by \cite[Lemma 3.5]{DS}, $\cS_\cZ$ is a tensor-suspended subcategory. Then the bijection trivially holds.
\end{proof}
\end{lem}

\begin{thm}\label{Aislebyaisle}
    For any commutative noetherian ring $R$ and finite acyclic quiver $Q$, it holds that
    \begin{enumerate}
        \item $\mathrm{Aisle_{cg}}_\boxtimes(\cD(RQ))\cong\Aisle_\mathrm{cg}(\cD(R))\times Q_0$;
        \item $\mathrm{Aisle_{cg}}_\boxtimes(\cD(RQ))\cong\Hom_{\Pos}(\Spec(R),\Aisle_\boxtimes(\cD(\bK Q)))$.
    \end{enumerate}
\begin{proof}
    (1) By \cite[Theorem 15 (i)]{DSw}, it suffices to prove that the following bijection holds
    \begin{equation*}\begin{gathered}
        \Susp_\boxtimes(\cD^c(RQ))\longleftrightarrow \Susp(\cD^c(R))\times Q_0 \\
        f:\cS\longmapsto\left(i^\ast\cS\right)_{i\in Q_0} \\
        \susp_{RQ}^\amalg\leftangle i_\times\cA_i\mid i\in Q_0\rightangle\longmapsfrom\left(\cA_i\right)_{i\in Q_0}:g
    \end{gathered}\end{equation*}
    Since both functors $i^\ast$ and $i_\times$ preserve compacts, to show that the assignments are well defined it is sufficient to show that, for any tensor-suspended subcategory $\cS$ in $\cD^c(RQ)$, $i^\ast\cS$ is extensions closed. Let $X, Z\in\cS$ and consider a distinguished triangle $i^\ast X\to M\to i^\ast Z\xrightarrow{+}\ $ in $\cD(R)$. In particular, applying the functor $i_\times$, by \cref{iprop}, we get that $X(i)\to i_\times M\to Z(i)\xrightarrow{+}\ $ is a distinguished triangle in $\cD(RQ)$ and, by tensor closure of $\cS$, $X(i), Z(i)\in\cS$. It follows that $i_\times M$ lies in $\cS$ and thus $M\cong i^\ast i_\times M\in i^\ast\cS$.

    Moreover, for a tensor-suspended subcategory $\cS\subseteq\cD^c(RQ)$, we have the equality
    \[g\circ f(\cS)=\susp_{RQ}^\amalg\leftangle i_\times i^\ast\cS\mid i\in Q_0\rightangle\]
    which, by \cref{iprop}, is equal to $\susp_{RQ}^\amalg\leftangle U(i)\boxtimes_{RQ}\cS\mid i\in Q_0\rightangle$ and, by \cref{AisleFilt} (1), it is equal to $\cS$. While, for a collection of tensor-suspended subcategories $\left(\cA_i\right)_{i\in Q_0}$ in $\cD(R)$, we have
    \[f\circ g(\left(\cA_i\right)_{i\in Q_0})=\left(i^\ast\susp_{RQ}^\amalg\leftangle j_\times\cA_j\mid j\in Q_0\rightangle\right)_{i\in Q_0}\]
    which, since we proved that $i^\ast$ reflects extensions, is equal to $\left(\susp_R^\amalg\leftangle i^\ast j_\times\cA_j\mid j\in Q_0\rightangle\right)_{i\in Q_0}$ and, by \cref{iprop}, it is equal to $(\cA_i)_{i\in Q_0}$.

    \noindent(2) Denote by $L(i)$ the $\bK Q$-simple module at a vertex $i\in Q_0$. By \cite[Theorem 15 (i)]{DSw}, it suffices to prove that the following bijection holds
    \begin{equation*}\begin{gathered}
        \Susp_\boxtimes(\cD^c(RQ))\longleftrightarrow \Hom_{\Pos}(\Spec(R),\Susp_\boxtimes(\cD^c(\bK Q))) \\
        f:\cS\longmapsto\left(\p\mapsto\susp_{\bK Q}^\amalg\leftangle  L(i)[n]\mid i_\times K(\p)[n]\in\cS\rightangle\right) \\
        \susp_{RQ}^\amalg\leftangle i_\times K(\p)[n]\mid L(i)[n]\in\sigma(\p),\p\in\Spec(R)\rightangle\longmapsfrom\,\sigma:g
    \end{gathered}\end{equation*}
    By \cref{ThickFilt} and \cref{OverField}, both $f(\cU)(\p)$ and $g(\sigma)$ are tensor-suspended subcategories. Moreover by \cite[Lemma 3.7]{Sab}, for any prime ideals $\q\subseteq\p$, we have that $K(\p)\in\susp_R^\amalg\leftangle K(\q)\rightangle$, thus $i_\times K(\p)[n]\in\susp_{RQ}^\amalg\leftangle i_\times K(\q)[n]\rightangle$ and $f(\cU)$ is a poset homomorphism. Since, by \cref{ThickFilt}, any tensor-suspended subcategory of $\cD^c(RQ)$ is generated by complexes of the form $i_\times K(\p)[n]$, we have that $g\circ f(\cS)=\cS$. Analogously, in virtue of the bijection between Serre subcategories and sets of simple modules (\cite[Section 8]{Kan}), by \cref{OverField}, we have that any tensor-suspended subcategory of $\cD^c(\bK Q)$ is generated by complexes of the form $L(i)[n]$ and so $f\circ g(\sigma)(\p)=\sigma(\p)$ for any $\p\in\Spec(R)$.
\end{proof}
\end{thm}

\begin{rem}
    Comparing \cref{Aislebyaisle} (2) with the classification in \cite{Sab}, where the author proves that $\Aisle_\mathrm{cg}(\cD(RQ))\cong\Hom_{\Pos}(\Spec(R),\Filt(\mathbf{Nc}(Q))$, one can note that the gap between compactly generated tensor-aisles and compactly generated aisles in $\cD(RQ)$ is reflected by the gap between the lattices $\Serre(\bK Q)\cong\cP(Q_0)$ (\cite[Section 8]{Kan}) and $\Wide(\bK Q)\cong\mathbf{Nc}(Q)$ (\cite[Theorem 1.1]{IT}).
\end{rem}

Now we get the analogous of \cref{Thickbyspc}, for compactly generated tensor-aisles.

\begin{cor}\label{Aislebyspc}
    For any commutative noetherian ring $R$ and finite acyclic quiver $Q$, we have that
    \[\mathrm{Aisle_{cg}}_\boxtimes(\cD(RQ))\cong\Filt(\cV(\Spc(\cD^c(RQ))))\]
    where $\cV(\Spc(\cD^c(RQ)))$ is the set of specialization closed subsets of $\Spc(\cD^c(RQ))$.
\begin{proof}
    Denote by $\cV(\Spec(R))$ the set of specialization closed subsets of the prime spectrum of $R$. By \cite{AJS}, the lattice of compactly generated aisles of $\cD(R)$ is isomorphic to $\Filt(\cV(\Spec(R)))$. Thus, since $Q_0$ is a discrete poset, by \cref{Aislebyaisle} (1), the lattice $\mathrm{Aisle_{cg}}_\boxtimes(\cD(RQ))$ is isomorphic to 
    \[\Hom_{\Pos}(Q_0,\Hom_{\Pos}(\bZ,\Hom_{\Pos}(\Spec(R),\{0,1\})))\]
    while, by \cref{Spectrum}, the lattice $\Filt(\cV(\Spc(\cD^c(RQ))))$ is isomorphic to
    \[\Hom_{\Pos}(\bZ,\Hom_{\Pos}(\Spec(R)\times Q_0,\{0,1\})))\]
    Note that applying the standard adjunction in $\Pos$ we can pass from one lattice to the other.
\end{proof}
\end{cor}

We want to end this section with the following open question arising from \cref{Aislebyspc}.
\begin{qst}
    Given a compactly generated tt-category $(\cT, \otimes, \mathbf{1})$ together with a ``standard'' aisle $\cT^{\leq 0}$ (the appropriate notion of \emph{standard} being part of the question), when are the compactly generated tensor-t-structures parametrized by filtrations of Thomason subsets of the Balmer spectrum ?
\end{qst}
Note that \cite{Hrb} and \cite{DS} answer affirmatively to this question for derived categories of commutative rings and noetherian schemes, respectively.

\section{Tensor telescope conjecture}\label{Sec4}
In the general setting of representations of an arbitrary quiver over commutative noetherian rings we do not know very much about compactly generated and homotopically smashing t-structures and how they interact. However, restricting to the tensor-compatible ones allows a classification of the former. We will see that this restriction is enough also for proving that any homotopically smashing tensor-t-structure is compactly generated.

\begin{lem}\label{SetGen}
    For an arbitrary ring $A$, any homotopically smashing t-structure $(\cU,\cV)$ in $\cD(A)$ is generated by objects, i.e.~there exists a set of objects $\cX\subseteq\cD(A)$ such that $\cU=\susp_{A}^\amalg\langle\cX\rangle$.
\begin{proof}
    Let us note first that, following the proof of \cite[Theorem 3.2 (1)]{MZ}, the homotopically smashing condition implies that the truncation functor $\tau_\cU^<$ of the t-structure $(\cU,\cV)$ commutes with directed homotopy colimits. Moreover, by the existence of semi-projective resolutions and \cite[Theorem 1.1]{CH}, any complex of $A$-modules is quasi-isomorphic to a semi-projective complex which is isomorphic to a directed homotopy colimit of compact objects $\{S_i\}_{i\in I}$ in $\cD(A)$. Then, for any $U\in\cU$, we have that
    \[U\cong\tau_\cU^<(U)\cong\tau_\cU^<(\hocolim_{i\in I}S_i)\cong\hocolim_{i\in I}\tau_\cU^<(S_i)\]
    In particular, it follows that $\cU=\susp_{RQ}^\amalg\langle\tau_\cU^<(S)\mid S\in\cD^c(RQ)\rangle$.
\end{proof}
\end{lem}

Recall that, by \cite[Theorem A]{SSV}, we can lift any t-structure $(\cM,\cN)$ in $\cD(R)$ to a t-structure $(\cM_Q,\cN_Q)$ in $\cD(RQ)$ defined as follows
\[\cM_Q=\{X\in\cD(RQ)\mid j^\ast X\in\cM \text{ for any } j\in Q\}\]
\[\cN_Q=\{X\in\cD(RQ)\mid j^\ast X\in\cN \text{ for any } j\in Q\}\]
note that the lifting is compatible with the orthogonality, in the sense that $(\cM_Q)^\perp=(\cM^\perp)_Q$.

In the following proposition, we show how a t-structure $(\cU,\cV)$ in $\cD(RQ)$ can be restricted to a t-structure in $\cD(R)$ preserving the homotopically smashing property. Specifically, the restriction is obtained by restricting the aisle and then taking its right orthogonal in $\cD(R)$. As illustrated in \cref{NaiveExmp}, the naive attempt to restrict both the aisle and the coaisle fails to even produce a t-structure. Therefore, this construction does not yield an inverse to the lifting in \cite{SSV}.

\begin{prop}\label{HomSm}
    For any homotopically smashing t-structure $(\cU,\cV)$ of $\cD(RQ)$, the pair $\left(i^\ast\cU,(i^\ast\cU)^\perp\right)$ is an homotopically smashing t-structure in $\cD(R)$. Moreover, its associated left truncation functor is given by the composite $i^\ast\circ\tau_Q^<\circ i_\ast$, where $\tau_Q^<$ is the truncation functor associated to the aisle $(i^\ast\cU)_Q$.
\begin{proof}
    Let $(\cU,\cV)$ a homotopically smashing t-structure in $\cD(RQ)$. By \cref{SetGen}, there exists a set of objects $\cX\subseteq\cD(RQ)$ such that $\cU=\susp_{RQ}^\amalg\langle X\mid X\in\cX\rangle$. By \cite[Theorem 3.4]{AJSo}, it follows that $i^\ast\cU=\susp_R^\amalg\langle i^\ast X\mid X\in\cX,i\in Q_0\rangle$ is an aisle in $\cD(R)$. We see that $(i^\ast\cU)^\perp$ is closed under directed homotopy colimits. Indeed, let $\{Y_\lambda\}_{\lambda\in\Lambda}\in\cD(R)^\Lambda$ be a coherent diagram such that $Y_\lambda\in(i^\ast\cU)^\perp$ for any $\lambda\in\Lambda$. In particular, we have that
    \[\Hom_{\cD(RQ)}(\cU,i_\ast Y_\lambda)\cong\Hom_{\cD(R)}(i^\ast\cU,Y_\lambda)=0 \text{ for any } \lambda\in\Lambda\]
    Since $\cU^\perp$ is closed under taking directed homotopy colimits and $i_\ast$ commute with them (see \cref{DerSec}), it follows that
    \[\Hom_{\cD(R)}(i^\ast\cU,\hocolim_{\lambda\in\Lambda}Y_\lambda)\cong\Hom_{\cD(RQ)}(\cU,\hocolim_{\lambda\in\Lambda}i_\ast Y_\lambda)=0\]
    Thus $\left(i^\ast\cU,(i^\ast\cU)^\perp\right)$ is an homotopically smashing t-structure in $\cD(R)$. As for the last part, let $U\in i^\ast\cU$ and $V\in\cD(R)$, then we have that $\Hom_{\cD(R)}(U,i^\ast\circ\tau_Q^<\circ i_\ast V)\cong\Hom_{\cD(RQ)}(i_!U,\tau_Q^<\circ i_\ast V)$. Since, by \cref{GroRem}, $i_!U$ lies in $(i^\ast\cU)_Q$ and $i^\ast i_!U\cong U$, the above is further equal to
    \[\Hom_{\cD(RQ)}(i_!U,i_\ast V)\cong\Hom_{\cD(R)}(i^\ast i_!U,V)\cong\Hom_{\cD(R)}(U,V)\]
    Thus, it is clear that the composite $i^\ast\circ\tau_Q^<\circ i_\ast$ is right adjoint to the inclusion $i^\ast\cU\hookrightarrow\cD(R)$.
\end{proof}
\end{prop}

\begin{exmp}\label{NaiveExmp}
    Consider the algebra $\bZ A_2$ and the t-structure $(\cU,\cV)$ generated by the representation $\bZ\xrightarrow{\pi}\bZ/p\bZ$ for some prime $p\in\bZ$. Note that, since there are no non-zero morphisms from $\bZ/p\bZ$ to $\bZ$, the representation $\bZ\xrightarrow{1}\bZ$ lies in $\cV$. In particular, the free module $\bZ$ is both in $1^\ast\cU$ and in $1^\ast\cV$, thus $(1^\ast\cU,1^\ast\cV)$ cannot be a t-structure.
\end{exmp}

We can now prove the tensor telescope conjecture for the tt-category $(\cD(RQ),\boxtimes_{RQ},\mathbf{U})$.

\begin{thm}\label{TTCthm}
    Any homotopically smashing tensor-t-structure of $\cD(RQ)$ is compactly generated.
\begin{proof}
    Let $(\cU,\cV)$ a homotopically smashing tensor-t-structure in $\cD(RQ)$. By \cref{HomSm}, we have that $\left(i^\ast\cU,(i^\ast\cU)^\perp\right)$ is an homotopically smashing t-structure in $\cD(R)$, thus by the telescope conjecture in $\cD(R)$ (see \cite{HN}), we can assume that $i^\ast\cU$ is a compactly generated aisle for any $i\in Q_0$. By the bijection \cref{Aislebyaisle} (1), $\susp_{RQ}^\amalg\langle i_\times i^\ast\cU\mid i\in Q_0\rangle$ is a compactly generated aisle in $\cD(RQ)$. By \cref{iprop}, we have that the above is equal to $\susp_{RQ}^\amalg\langle X(i)\mid X\in\cU, i\in Q_0\rangle$ and, by \cref{AisleFilt}, it is further equal to $\cU$.
\end{proof}
\end{thm}

\subsection{Other tensor-t-structures}
Actually, we can strengthen the tensor telescope conjecture of \cref{TTCthm}, by generalizing the result to homotopically smashing tensor-t-structures with respect to a broad class of suspended subcategories $\cT^{\leq0}$, which are strictly contained in the standard aisle.

For an $RQ$-module $C$, the \emph{support of $C$ over $Q$} is the full subquiver of $Q$ with vertices $i\in Q_0$ such that $C_i$ is non-zero and we denote it by $\Supp_Q(C)$.

\begin{defn}\label{FiltSyst}
We define a \emph{filtration system of the unit} to be a set of $R$-free $RQ$-modules $\cC=\{C_1,\ldots,C_n\}$ such that:
\begin{itemize}
    \item[(i)\,] $\mathbf{U}\in\mathsf{filt}\langle\cC\rangle$;
    \item[(ii)] $\Supp_Q(C_i)\cap\Supp_Q(C_j)=\emptyset$ for any $i\neq j$.
\end{itemize}
Moreover, we say that a filtration system (of the unit) has \emph{Dynkin support} if, for any $k=1,\ldots,n$, the quiver $Q_k:=\Supp_Q(C_k)$ is a Dynkin quiver.
\end{defn}

Given a filtration system $\cC$, we say that a suspended subcategory $\cS$ is \emph{$\cC$-tensor-suspended} if it is tensor-suspended with respect to $\susp_{RQ}\langle\cC\rangle$, i.e.~by \cite[Lemma 3.5]{DS} if $\cC\boxtimes_{RQ}\cS\subseteq\cS$. A t-structure is a \emph{$\cC$-tensor-t-structure} if the aisle is a $\cC$-tensor-suspended subcategory. Finally, given a set of objects $\cX\subseteq\cT$, we denote the smallest cocomplete $\cC$-tensor-suspended subcategory containing $\cX$ by $\langle\cX\rangle^\cC$, by $\Susp_\boxtimes^\cC(\cD^c(RQ))$ the lattice of $\cC$-tensor-suspended subcategories of $\cD^c(RQ)$ and by $\mathrm{Aisle_{cg}^\cC}_\boxtimes(\cD(RQ))$ the lattice of compactly generated $\cC$-tensor-aisle of $\cD(RQ)$.
\pagebreak

\begin{exmp}\label{DynkinExmp}\hfill
    \begin{enumerate}
        \item Any t-structure in $\cD(RQ)$ is a $\cC$-tensor-t-structure with respect to the trivial filtration system $\cC=\{\mathbf{U}\}$. In particular, when $Q$ is a Dynkin quiver, any t-structure is a $\cC$-tensor-t-structure with respect to a filtration system with Dynkin support.
        \item The tensor-t-structures with respect to the standard aisle $\cD^{\leq0}(RQ)$ are precisely the $\cC$-tensor-t-structures with respect to the filtration system $\cC=\{U(i)\mid i\in Q_0\}$, which always has Dynkin support.
    \end{enumerate} 
\end{exmp}

\begin{rem}
    Given a filtration system $\cC$, any module $C_k\in\cC$ defines a fully faithful functor $c_k:Q_k\hookrightarrow Q$ and thus an adjoint triple ${c_k}_!\dashv{c_k}^\ast\dashv{c_k}_\ast$ as follows
    \[\begin{tikzcd} \cD(RQ) \arrow[rr, "{c_k}^\ast"] & & \cD(RQ_k) \arrow[ll, "{c_k}_!"', bend right] \arrow[ll, "{c_k}_\ast", bend left] \end{tikzcd}\]
    Since $c_k$ is fully faithful, by the characterization of the projective $RQ$-modules in \cite[Theorem 3.1]{EE}, given a bounded complexes $P$ of finitely generated projective $RQ$-modules, ${c_k}^\ast P$ is a bounded complexes of finitely generated projective $RQ_k$-modules. Thus, ${c_k}^\ast$ preserves compact objects and, since coproducts in $\cD(RQ)$ are computed vertexwise, it also preserves coproducts. In particular, as in \cref{GroRem}, ${c_k}_!$ preserves compacts, ${c_k}_\ast$ preserves coproducts and all three ${c_k}_!, {c_k}^\ast, {c_k}_\ast$ preserve directed homotopy colimits.
    Moreover, by \cite[Proposition 1.26]{Gro}, the composite ${c_k}^\ast{c_k}_!$ is natural isomorphic to $\id_{\cD(RQ_k)}$. So, we can define the functor
    \[{c_k}_\times:=C_k\boxtimes_{RQ}{c_k}_!:\cD(RQ_k)\to\cD(RQ)\]
\end{rem}

In view of \cref{DynkinExmp}, the following proposition generalizes \cref{AisleFilt} (1) and \cref{Aislebyaisle} (1) to $\cC$-tensor-t-structures.

\begin{prop}\label{Aislebyfilt}
    For any filtration system of the unit $\cC=\{C_1,\ldots,C_n\}$, it holds that:
    \begin{enumerate}
        \item $\langle X\rangle^\cC=\susp_{RQ}^\amalg\langle C_k\boxtimes_{RQ}X\mid k=1,\ldots,n\rangle$;
        \item $\mathrm{Aisle_{cg}^\cC}_\boxtimes(\cD(RQ))\cong\prod\limits_{k=1}^n\Aisle_\mathrm{cg}(\cD(RQ_k))$
    \end{enumerate}
\begin{proof}
    (1) Since the tensor unit $\mathbf{U}$ admits a filtration whit factors isomorphic to the representations $C_k$, tensoring this filtration by a complex $X\in\cD(RQ)$ it follows that $X\in\susp_{RQ}^\amalg\leftangle C_k\boxtimes_{RQ}X\mid k=1,\ldots,n\rightangle$. Moreover, the latter is a $\cC$-tensor-suspended subcategory by \cite[Lemma 3.5]{DS}. On the other hand, by closure condition, also the reverse inclusion holds.

    \noindent(2) By \cite[Theorem 15 (i)]{DSw}, it suffices to prove that the following bijection holds
    \begin{equation*}\begin{gathered}
        \Susp_\boxtimes(\cD^c(RQ))\longleftrightarrow \prod\limits_{k=1}^n\Susp(\cD^c(RQ_k))\\
        f:\cS\longmapsto \left({c_k}^\ast\cS\right)_{k=1,\ldots,n} \\
        \susp_{RQ}^\amalg\leftangle{c_k}_\times\cA_k\mid k=1,\ldots,n\rightangle\longmapsfrom\left(\cA_k\right)_{k=1,\ldots,n}:g
    \end{gathered}\end{equation*}
    Since both functors ${c_k}^\ast$ and ${c_k}_\times$ preserve compacts and ${c_k}^\ast$ preserves extensions, shifts and summands, the assignments are well defined. Moreover, for a tensor-suspended subcategory $\cS\subseteq\cD(RQ)$, the composite $g\circ f(\cS)$ is $\susp_{RQ}^\amalg\leftangle{c_k}_\times{c_k}^\ast\cS\mid k=1,\ldots,n\rightangle$. By \cref{iprop}, we have that 
    \[g\circ f(\cS)=\susp_{RQ}^\amalg\leftangle C_k\boxtimes_{RQ}\cS\mid k=1,\ldots,n\rightangle\]
    which, by point (1), is equal to $\cS$. For a family of tensor-ideals $\left(\cA_k\right)_{k=1,\ldots,n}$ with $\cA_k\subseteq\cD(RQ_k)$, the composite $f\circ g(\left(\cA_k\right)_{k=1,\ldots,n})$ is given by $\left({c_k}^\ast\susp_{RQ}^\amalg\leftangle {c_\ell}_\times\cA_\ell\mid \ell=1,\ldots,n\rightangle\right)_{k=1,\ldots,n}$ which, by the properties of ${c_k}^\ast$ and the fact that the supports of objects in $\cC$ are disjoint, is equal to $\left(\susp_R^\amalg\leftangle {c_k}^\ast{c_\ell}_\times\cA_\ell\mid\ell=1,\ldots,n\rightangle\right)_{k=1,\ldots,n}=(A_k)_{k=1,\ldots,n}$.
\end{proof}
\end{prop}

We can now prove the tensor telescope conjecture for $\cC$-tensor-t-structures with respect to some filtration system with Dynkin support and so, generalize \cref{TTCthm}.

\begin{thm}\label{TTCimprove}
    Given a filtration system with Dynkin support $\cC$, any homotopically smashing $\cC$-tensor-t-structure of $\cD(RQ)$ is compactly generated.
\begin{proof}
    Let $(\cU,\cV)$ a homotopically smashing $\cC$-tensor-t-structure in $\cD(RQ)$. By similar arguments to those in \cref{HomSm}, we have that $\left({c_k}^\ast\cU,({c_k}^\ast\cU)^\perp\right)$ is a homotopically smashing t-structure in $\cD(RQ_k)$ for any $k=1,\ldots,n$, thus by the telescope conjecture in $\cD(RQ_k)$ (see \cite{Sab}), we can assume that ${c_k}^\ast\cU$ is a compactly generated aisle for any $k=1,\ldots,n$. By the bijection in \cref{Aislebyfilt} (2), $\susp_{RQ}^\amalg\langle{c_k}_\times{c_k}^\ast\cU\mid k=1,\ldots,n\rangle$ is a compactly generated $\cC$-tensor-aisle in $\cD(RQ)$. Moreover, it is equal to $\susp_{RQ}^\amalg\langle C_k\boxtimes_{RQ}\cU\mid k=1,\ldots,n\rangle$ and thus, by \cref{Aislebyfilt} (2), to $\cU$.
\end{proof}
\end{thm}

The next example illustrates how \cref{TTCimprove} could represent a significant step toward proving the telescope conjecture for representations of finite, acyclic, and simply-laced quivers.

\begin{exmp}
    Consider the Euclidean quiver
    \begin{center}
        $\widetilde{D}_5=$\scriptsize
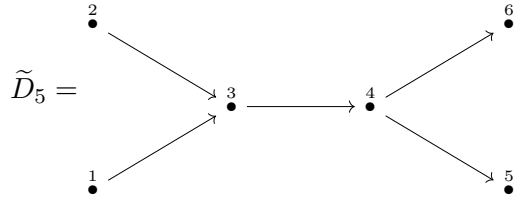
        \begin{tikzcd}[ampersand replacement=\&]
            \overset{2}{\bullet} \arrow[rrd] \& \& \& \& \& \& \overset{6}{\bullet} \\
            \& \& \overset{3}{\bullet} \arrow[rr] \& \& \overset{4}{\bullet} \arrow[rru] \arrow[rrd] \& \& \\
            \overset{1}{\bullet} \arrow[rru] \& \& \& \& \& \& \overset{5}{\bullet}
        \end{tikzcd}
    \end{center}

    \noindent Let $D_5$ be the full subquiver spanned by the vertices $\{1,2,3,4,5\}$ and $P$ the representation given by $P_i=R$ and $P_\alpha=\id_R$ if $i,\alpha\in D_5$ and $P_6=0$, then $\cC=\{P,U(6)\}$ is a filtration system with Dynkin support. Note that a suspended subcategory $\cU$ is a $\cC$-tensor-ideal if and only if for any object $X\in\cU$ such that $\Supp_Q(X)\nsubseteq D_5$ then $X(6)\in\cU$. Thus, by \cref{TTCimprove}, any homotopically smashing t-structure whose aisle satisfies this property is compactly generated.

    Moreover, in this example as well as in any Euclidean quiver, a non-trivial filtration systems $\cC\neq\{\mathbf{U}\}$ has always Dynkin support. Therefore, for a homotopically smashing t-structure, it is enough to have a $\cC$-tensor-aisle for some non-trivial filtration system $\cC$ in order to be compactly generated.
\end{exmp}

\raggedbottom
\bibliographystyle{amsalpha}
\bibliography{references}
\end{document}